\documentclass[12pt,a4paper]{article}
\usepackage{latexsym}
\usepackage{mathrsfs}
\usepackage{amsmath}
\usepackage{amssymb}
\usepackage[usenames,dvipsnames,svgnames,table]{xcolor}

\newtheorem{theorem}{Theorem}[section]
\newtheorem{lemma}[theorem]{Lemma}

\def \H {{\cal{H}}}

\def \proof {\noindent{\bf Proof}\quad}

\newcommand{\qed}{\hfill$\Box$\vspace{0.2cm}}

\title{On Hamilton Decompositions of Line Graphs\\ of Non-Hamiltonian Graphs and\\ Graphs without Separating Transitions}

\author{
Darryn Bryant, Barbara Maenhaut, Benjamin R. Smith
}

\date{ }

\begin{document}
\maketitle\thispagestyle{empty}
\def\baselinestretch{1.25}\small\normalsize

\vspace{-1cm}

Department of Mathematics, University of Queensland, QLD 4072, Australia.\newline 
Email: \texttt{db@maths.uq.edu.au}, \quad \texttt{bmm@maths.uq.edu.au},\quad  \texttt{bsmith.maths@gmail.com}

\vspace{0.5cm}

\begin{abstract}
In contrast with Kotzig's result that the line graph of a $3$-regular graph $X$ is Hamilton decomposable if and only if $X$ is Hamiltonian,
we show that for each integer $k\geq 4$ there exists a simple non-Hamiltonian $k$-regular graph whose line graph has a Hamilton decomposition.
We also answer a question of Jackson by showing that for each integer $k\geq 3$ there exists a simple connected $k$-regular graph with no separating transitions whose line graph has no Hamilton decomposition. 
\end{abstract}

\section{Introduction}

In the 1960's Kotzig \cite{Kot} proved that the existence of a Hamilton cycle in a $3$-regular graph $X$ is both necessary and sufficient for the existence of a Hamilton decomposition of its line graph $L(X)$. Hamilton decomposability of line graphs has subsequently been studied extensively, but the general question of classifying those graphs whose line graphs have Hamilton decompositions remains open. 
This topic has been considered from a number of different perspectives. In particular, Hamilton decomposability of $L(X)$ has been considered with imposed conditions on the connectivity \cite{FleHilJac,Jac,JacWor} or Hamiltonicity \cite{Ber,Jae,MutPau,Pik1,Pik2,Pik3} of $X$. 
Additional papers containing results related to Hamilton decompositions of line graphs include \cite{HeiVer,Ver} and the survey \cite{AlsBerSot}. 

In this paper we answer a question of Jackson \cite{Jac} on Hamilton decomposability of the line graphs of graphs with no separating transitions (a connectivity-related condition defined below), 
and we prove that the above-mentioned result of Kotzig does not hold for $k$-regular graphs when $k\geq 4$.
If $X$ is regular of degree $2k$ or $2k+1$, then a set of $k$ pairwise edge-disjoint Hamilton cycles in $X$ is 
called a {\em Hamilton decomposition}, and a graph admitting a Hamilton decomposition is said to be {\em Hamilton decomposable}. 

In \cite{Jac}, Jackson calls a pair of half edges incident with a vertex $u$ a {\em transition} at $u$, and 
if $t$ is a transition at $u$
in a graph $X$, then he defines $X^t$ to be the graph obtained from $X$ by splitting $u$ into two new vertices $u_1$ and $u_2$, 
joining the two half edges of $t$ to $u_1$, and joining each other half edge at $u$ to $u_2$.
A {\em separating transition} is then defined to be a transition $t$ such that $X^t$ has more components than $X$. 
It is shown in \cite{Jac} that for $k\geq 3$, 
the line graph of a connected $k$-regular graph $X$ is $(2k-2)$-edge-connected if and only if $X$ has no 
separating transitions (the result is actually stated only for the case $k$ is even, but the same argument works when $k$ is odd).
It follows that if $k\geq 3$ and $X$ is any $k$-regular graph with a separating transition, then $L(X)$ has no Hamilton decomposition.  
We observe that the preceding statement is not true without the requirement that $X$ be regular. For example, any star with at least three edges has both a separating transition and a Hamilton decomposable line graph.

Having observed that absence of separating transitions in $X$ is necessary for Hamilton decomposability of $L(X)$, 
Jackson asks (Problem 5.2 in \cite{Jac}) whether it is true that the line graph of a connected $2k$-regular graph $X$ has a Hamilton decomposition if and only if $X$ has no separating transitions. 
The same question could also be asked for connected regular graphs of odd degree. 
The answer is no in the case of $3$-regular graphs because 
there are many connected non-Hamiltonian $3$-regular graphs that have no separating transitions, 
and Kotzig's result tells us the the line graphs of these graphs have no Hamilton decomposition.  
In Section \ref{septrans} we construct for each integer $k\geq 3$ 
a simple connected $k$-regular graph with no separating transitions whose line graph has no Hamilton decomposition,
thereby showing that the answer to Jackson's question is no for every degree greater than $3$.

The authors \cite{BryMaeSmi} have recently shown that the existence of a Hamilton cycle in a simple graph 
$X$ is sufficient for Hamilton decomposability of $L(X)$ when $X$ is regular of even degree, and that the existence of a Hamiltonian $3$-factor in $X$ is 
sufficient for Hamilton decomposability of $L(X)$ when $X$ is regular of odd degree. 
Whether the existence of a Hamilton cycle,
rather than a Hamiltonian $3$-factor, is sufficient for Hamilton decomposability of $L(X)$ when $X$ is regular of odd degree remains
an open question. 
The results just mentioned partially extend Kotzig's result to $k$-regular graphs with $k>3$, but only in the direction of sufficiency.
Going in the opposite direction, we show in Section \ref{nonHam} that the existence of a Hamilton cycle in $X$ is not necessary for Hamilton decomposability of $L(X)$ when $X$ is regular of degree at least $4$.

The proofs of both of our main results involve construction of new graphs by deletion of an edge of a graph and insertion of the resulting graph into an edge of another graph, and we now give the formal definition of this procedure. 
Let $X$ and $X'$ be vertex-disjoint graphs (not necessarily simple), let $u$ and $v$ be adjacent vertices in $X$, and let $u'$ and $v'$ be adjacent vertices in $X'$.
We define the {\em insertion} of $X'-u'v'$ into an edge $uv$ of $X$ to be the graph obtained 
from $X\cup X'$ by replacing an edge $uv$ of $X$ and an edge $u'v'$ of $X'$ with an edge joining $u$ to $u'$ and an edge joining $v$ to $v'$. 
In this definition the order in which the vertices of the edges $uv$ and $u'v'$ are listed may change the resulting graph, 
but this will be of no consequence in our constructions. 

\section{Separating transition-free graphs whose\\ line graphs are not Hamilton decomposable}\label{septrans}

\begin{theorem}\label{septransthm}
For each integer $k\geq 3$, there exists a simple connected $k$-regular graph with no separating transitions whose line graph has no Hamilton decomposition.
\end{theorem}

\proof
For each integer $k\geq 3$ and each even integer $t\geq 4$, define $Y_{k,t}$ to be the multigraph with vertices $v_1,v_2,\ldots,v_t$, 
and edge set given by joining $v_i$ to $v_{i+1}$ with two edges for $i=1,3,\ldots,t-1$, and joining $v_i$ to $v_{i+1}$ 
with $k-2$ edges for $i=2,4,\ldots,t$. Here, and throughout what follows, $v_{t+1}$ is identified with $v_{1}$. 
Let $X_{k,t}$ be the graph obtained from $Y_{k,t}$ by inserting a copy of $K_{k+1}-e$ into each edge of $Y_{k,t}$.
It is easy to see that $X_{k,t}$ is a simple $k$-regular graph that has no separating transitions. We will 
show that $L(X_{k,t})$ has no Hamilton decomposition for $t\geq k$, but first we need to introduce labels for various edges 
of $X_{k,t}$.

For $i=1,3,\ldots,t-1$, let $X_i^1$ and $X_i^2$ be the two copies of $K_{k+1}-e$ that are inserted into the two edges joining $v_i$ to $v_{i+1}$. 
For $i=1,2,\ldots,t$, let $e_i^1,e_i^2,\ldots,e_i^k$ be the $k$ edges of $X_{k,t}$ that are incident with $v_i$.
For $i=1,3,\ldots,t-1$ and for $j=1,2$, let $e_i^j$ be the unique edge joining $v_i$ to $X_i^j$,
and let $e_{i+1}^j$ be the unique edge joining $v_{i+1}$ to $X_i^j$.
For $i=1,3,\ldots,t-1$ and for $j=1,2$, 
let $f_{i,1}^j,f_{i,2}^j,\ldots,f_{i,k-1}^j$ be the $k-1$ edges of $X_i^j$ that are adjacent to $e_i^j$,
and
let $f_{i+1,1}^j,f_{i+1,2}^j,\ldots,f_{i+1,k-1}^j$ be the $k-1$ edges of $X_i^j$ that are adjacent to $e_{i+1}^j$.
Finally, for $i=1,2,\ldots,t$, let $E_i$ be the set of $2(k-2)$ edges of $L(X_{k,t})$ having one endpoint in $\{e_i^1,e_i^2\}$ and the other
in $\{e_i^3,e_i^4,\ldots,e_i^d\}$.

For a contradiction, suppose $t\geq k$ and $\H$ is a Hamilton decomposition of $L(X_{k,t})$. 
Note that $\H$ contains $k-1$ Hamilton cycles. 
Since $t>k-1$, in $L(X_{k,t})$ at least two of the edges $e_1^1e_1^2,e_2^1e_2^2,\ldots,e_t^1e_t^2$ are in the same Hamilton cycle of $\H$.
Let this cycle be $H\in\H$ and let $e_a^1e_a^2$ and $e_b^1e_b^2$ be distinct edges of $H$ (so $a,b\in\{1,2,\ldots,t\}$). 

Now, for $i=1,3,\ldots,t-1$ and for $j=1,2$, $\{e_i^j,e_{i+1}^j\}$ is a vertex cut of $L(X_{k,t})$, and it follows that for $i=1,2,\ldots,t$ and for $j=1,2$, each Hamilton cycle of $\H$ contains exactly one of the $k-1$ edges $e_i^jf_{i,1}^j,e_i^jf_{i,2}^j,\ldots,e_i^jf_{i,k-1}^j$ of $L(X_{k,t})$. 
But this implies that $H$ contains none of the edges of $E_a$ and none of the edges of $E_b$. 
Since $E_a\cup E_b$ is an edge cut of $L(X_{k,t})$, this is a contradiction, and we conclude that $L(X_{k,t})$ has no Hamilton decomposition. 
\qed

\section{Non-Hamiltonian graphs whose line graphs are Hamilton decomposable}\label{nonHam}

Hamilton cycles in $L(X)$ are related to certain Euler tours of $X$.
If $$v_0,e_1,v_1,e_2,v_2,\ldots,e_t,v_t=v_0$$ is any Euler tour of $X$
(so the edge set of $X$ is $\{e_1,e_2,\ldots,e_t\}$ and $v_0,v_1,\ldots,v_t$ are vertices of $X$), then $(e_1,e_2,\ldots,e_k)$ is a Hamilton cycle in 
$L(X)$. However, not every Hamilton cycle in $L(X)$ corresponds to an Euler tour of $X$. For example, $(wz,wx,wy,yz,xy,xz,wz)$ is a Hamilton cycle in the line graph of the complete graph with vertex set $\{w,x,y,z\}$, but it clearly does not correspond to an Euler tour of this complete graph. Indeed, there is no Euler tour of the complete graph of order $4$.

We shall say that a Hamilton cycle in $L(X)$ is {\em Euler tour compatible} if it corresponds to an Euler tour in $X$.
In order to say more about the properties that Hamilton cycles in $L(X)$ must have in order that they be Euler tour compatible, we make the following definitions.  
If $u$ is a vertex in a graph $X$, then the neighbourhood of $u$ in $X$ is denoted by $N_X(u)$. 
Suppose $X$ is a simple graph and $N_X(u)=\{v,a_1,a_2,\ldots,a_k\}$. Then in $L(X)$ the {\em $u$-neighbourhood} of the vertex $uv$ is 
$\{ua_1,ua_2,\ldots,ua_k\}$ and is denoted by $N^u_{L(X)}(uv)$.
Thus, $N^u_{L(X)}(uv)\cap N^v_{L(X)}(uv)=\emptyset$ and $N_{L(X)}(uv)=N^u_{L(X)}(uv)\cup N^v_{L(X)}(uv)$.

It is easy to see that a Hamilton cycle $H$ in $L(X)$ is {\em Euler tour compatible} if and only if for each vertex $uv$ in $L(X)$,
one neighbour in $H$ of $uv$ is from the the $u$-neighbourhood of $uv$ and the other neighbour in $H$ of $uv$ is from the $v$-neighbourhood
of $uv$. If this propery holds for a vertex $uv$ in a Hamilton cycle $H$ of $L(X)$, 
then we say that $H$ is {\em Euler tour compatible at $uv$}. Thus, a Hamilton cycle is Euler tour compatible if and only if it is 
Euler tour compatible at each of its vertices. More generally, a Hamilton decomposition of $L(X)$ is Euler tour compatible at $uv$ if each of its Hamilton cycles is Euler tour compatible at $uv$, and is {\em everywhere Euler tour compatible} if it is Euler tour compatible at every vertex of $L(X)$.

A Hamilton decomposition of $L(X)$ that is everywhere Euler tour compatible is thus equivalent to a 
{\em perfect set of Euler tours} of $X$, 
where a set $S$ of Euler tours of $X$ is perfect if each 2-path in $X$ occurs in exactly one Euler tour in $S$. 
In \cite{HeiVer}, Heinrich and Verrall construct perfect sets of Euler tours for each complete graph of odd order, thus establishing the following theorem. 

\begin{theorem}\label{theoremA}[Heinrich and Verrall \cite{HeiVer}]
For each odd integer $n\geq 3$, the line graph of the complete graph of order $n$ has a Hamilton decomposition that is everywhere Euler tour compatible. 
\end{theorem}

The complete graph of even order has no Euler tour. However, there is a natural way to extend the above ideas by considering instead the multigraph $K_n+I$ which is obtained from the complete graph of even order $n$ by duplicating each edge in a set 
$I$ of edges that form a perfect matching. 
In \cite{Ver}, Verrall shows that $K_n+I$ has a perfect set of Euler tours for all even $n\geq 4$, where the definition of  
perfect set of Euler tours is suitably modified to accommodate the edges of multiplicity 2. 
The modification is exactly what is needed to ensure that perfect sets of Euler tours of $K_n+I$ correspond to 
Hamilton decompositions of $L(K_n)$ that are Euler tour compatible at each vertex of $L(K_n)$ except those in $I$.
Indeed, as stated in \cite{Ver}, the modification is made specifically to parallel Theorem \ref{theoremA}, and
it is easily verified that the main result 
in \cite{Ver} can be restated in our terminology as follows.  

\begin{theorem}\label{theoremB}[Verrall \cite{Ver}]
If $n\geq 4$ is an even integer, $K$ is a complete graph of order $n$, and $I$ is a perfect matching in $K$, then
$L(K)$ has a Hamilton decomposition that is Euler tour compatible at each vertex of $V(L(K))\setminus I$.
\end{theorem}

\begin{lemma}\label{lemmaA}
Let $X$ and $X'$ be vertex-disjoint $k$-regular graphs, let $uv$ be an edge in $X$, let $u'v'$ be an edge in $X'$, and let $Y$ be the insertion of $X'-u'v'$ into the edge $uv$ of $X$. 
If $\H$ is a Hamilton decomposition of $L(X)$ that is Euler tour compatible at $uv$ and $\H'$ is a Hamilton decomposition of $L(X')$ that is Euler tour compatible at $u'v'$, 
then there exists a Hamilton decomposition $\H^*$ of $L(Y)$ such that if $\H$ is Euler tour compatible at a vertex $xy\neq uv$ of $L(X)$, 
then $\H^*$ is also Euler tour compatible at $xy$.
\end{lemma}

\proof
Suppose $\H=\{H_1,H_2,\ldots,H_{k-1}\}$ is a Hamilton decomposition of $L(X)$ that is Euler tour compatible at $uv$ and suppose $\H'=\{H'_1,H'_2,\ldots,H'_{k-1}\}$ is a Hamilton decomposition of $L(X')$ that is Euler tour compatible at 
$u'v'$. 
For $i=1,2,\ldots,k-1$, let the two neighbouring vertices of $uv$ in $H_i$ be $ua_i$ and $vb_i$, let the 
two neighbouring vertices of $u'v'$ in $H'_i$ be $u'a'_i$ and $v'b'_i$, and let $J_i$ be the graph obtained from the union of $H_i$ and $H'_i$
by
replacing the vertices $uv$ and $u'v'$ with $uu'$ and $vv'$, 
replacing the edge joining $ua_i$ to $uv$ with an edge joining $ua_i$ to $uu'$, 
replacing the edge joining $vb_i$ to $uv$ with an edge joining $vb_i$ to $vv'$, 
replacing the edge joining $u'a'_i$ to $u'v'$ with an edge joining $u'a'_i$ to $uu'$, and 
replacing the edge joining $v'b'_i$ to $u'v'$ with an edge joining $v'b'_i$ to $vv'$.
It is easily seen that $\H^*=\{J_1,J_2,\ldots,J_{k-1}\}$ is the required Hamilton decomposition of $L(Y)$.
\qed

\begin{theorem}
For each integer $k\geq 4$ there exists a simple non-Hamiltonian $k$-regular graph whose line graph has a Hamilton decomposition. 
\end{theorem}

\proof
Let $k\geq 4$, let $v,u_1,u_2$ and $u_3$ be vertices in a complete graph $X$ of order $k+1$.
For each $i\in\{1,2,3\}$, let $X'_i$ be a complete graph of order $k+1$, let $u'_iv'_i$ be an edge in $X'_i$, 
and insert $X'_i-u'_iv'_i$ into the edge $u_iv$ of $X$.
Let $Y$ be the resulting graph. 
We claim that $Y$ is a non-Hamiltonian $k$-regular graph. To see that $Y$ is non-Hamiltonian, observe that for $i\in\{1,2,3\}$,
$\{vv'_i,u_iu'_i\}$ is
an edge cut, and so any Hamilton cycle necessarily contains the three edges $vv'_1$, $vv'_2$ and $vv'_3$, which is impossible. 

We now use Theorems \ref{theoremA} and \ref{theoremB} and Lemma \ref{lemmaA} to show that $L(Y)$ has a Hamilton decomposition. 
Since $k\geq 4$, $L(X)$ has a Hamilton decomposition that is Euler tour compatible at $vu_1$, $vu_2$ and $vu_3$ by Theorem \ref{theoremA} ($k$ even)
or \ref{theoremB} ($k$ odd). 
Also by Theorem \ref{theoremA} ($k$ even)
or \ref{theoremB} ($k$ odd), for each $i\in\{1,2,3\}$, $L(X'_i)$ has a Hamilton decomposition that is Euler tour compatible at $u'_iv'_i$.
It thus follows by Lemma \ref{lemmaA} (applied three times) that $L(Y)$ has a Hamilton decomposition. 
\qed

\vspace{0.5cm}

\noindent{\bf Acknowledgement:}
The authors acknowledge the support of the Australian Research Council via grants DP150100530, DP150100506,
DP120100790, DP120103067 and DP130102987.

\end{document}